\documentclass[reqno]{amsart}     
\usepackage[foot]{amsaddr}
\usepackage{amsmath,amssymb}    
\usepackage{enumerate}
\usepackage{bbm}
\usepackage[backend=bibtex,doi=false,url=false,firstinits=true, maxbibnames=10, defernumbers=true, style=numeric]{biblatex}
\usepackage{textcomp}

\newtheorem{theorem}                   {Theorem} 
\newtheorem{thm}             [theorem] {Theorem} 
\newtheorem{lemma}           [theorem] {Lemma}

\theoremstyle{remark}

\newcommand{\eps}{\varepsilon}
\newcommand{\cH}{\mathcal{H}}

\newcommand{\Real}{\mathbb{R}}

\newcommand{\conntime}{{\tau_c}}
\newcommand{\isoltime}{{\tau_i}}
\newcommand{\process}{{\{\cH^k(n,M)\}_M}}
\newcommand{\graphprocess}{{\{G(n,M)\}_M}}

\newcommand{\hknp}{{\cH^k(n,p)}}

\newcommand{\Po}{\mathrm{Po}}
\renewcommand{\Pr}{\mathbb{P}}

\newcommand{\giantthreshold}{p_g}
\newcommand{\connthreshold}{p_c}

\DeclareFieldFormat{pages}{#1}
    \DeclareFieldFormat[article]{date}{\mkbibparens{#1}}
     \DeclareFieldFormat[article]{volume}{\bf{#1}}
     \DeclareFieldFormat[article]{number}{\printtext{no}\adddot\addnbthinspace{#1}\addspace}

\renewbibmacro*{journal+issuetitle}{%
  \usebibmacro{journal}%
  \setunit*{\addnbspace}%
  \iffieldundef{series}
    {}
    {\newunit
     \printfield{series}%
     \setunit{\addspace}}%
  \usebibmacro{volume+date+number}%
   \newunit}

\newbibmacro*{volume+date}{%
    \iffieldundef{volume}
      {\usebibmacro{date}}
      {\printfield{volume}%
       \setunit*{\addthinspace}%
       \usebibmacro{date}}
 \newunit}

\newbibmacro*{volume+date+number}{%
    \iffieldundef{number}
      {\usebibmacro{volume+date}}
      {\usebibmacro{volume+date}}%
       \addcomma\addspace
       \printfield{number}%
  \newunit}

\AtBeginDocument{}

\renewbibmacro{in:}{%
 \ifentrytype{article}{}{\printtext{\bibstring{in}\intitlepunct}}}

 \AtEveryBibitem{\clearfield{issn}}
 \AtEveryCitekey{\clearfield{issn}}

\begin{document}
\title[Evolution of  random hypergraphs]{Evolution of  high-order connected components in random hypergraphs}
\thanks{The authors are supported by Austrian Science Fund (FWF): P26826, W1230.}
\thanks{\textrm{DOI:} http://dx.doi.org/10.1016/j.endm.2015.06.077. Copyright \textcopyright 2015. This manuscript version is made available under the CC-BY-NC-ND 4.0 license http://creativecommons.org/licenses/by-nc-nd/4.0/}

\author[O.~Cooley, M.~Kang and C.~Koch]{Oliver Cooley, Mihyun Kang and Christoph Koch}
\email{\{cooley,kang,ckoch\}@math.tugraz.at}
\address{Institute of Optimization and Discrete Mathematics, Graz University of Technology, 8010 Graz, Austria}

\date{}

\begin{abstract}
We consider high-order connectivity in $k$-uniform hypergraphs defined as follows: Two $j$-sets are $j$-connected if there is a walk of edges between them such that two consecutive edges intersect in at least $j$ vertices. We describe the evolution of $j$-connected components in the $k$-uniform binomial random hypergraph $\hknp$. In particular, we determine the asymptotic size of the giant component shortly after its emergence and establish the threshold at which the $\hknp$ becomes $j$-connected with high probability. We also obtain a hitting time result for the related random hypergraph process $\process$ -- the hypergraph becomes $j$-connected exactly at the moment when the last isolated $j$-set disappears. This generalises well-known results for graphs and vertex-connectivity in hypergraphs.
\end{abstract}

\maketitle

\section{Evolution of random graphs}

The theory of random graphs was founded in the late 1950s by Erd\H{o}s and R\'enyi describing the evolution of the random graph process $\graphprocess.$ The vertex set of this process is $[n]:=\left\{1,\dots,n\right\}$ and initially there are no edges present. In each step of the process, add an edge between a pair of vertices chosen uniformly at random amongst all pairs of vertices that do not already form an edge. In the early stages of this process, all connected components are small and then, within very short time, they merge into a single component of linear size -- the \emph{giant component}. This remarkable phenomenon, first proved in~\cite{ErdosRenyi60}, is known as the \emph{phase transition} of the random graph process $\graphprocess$.

It is often more convenient to analyse the binomial random graph $G(n,p)$: The vertex set is $[n]$ and every pair of vertices is connected by an edge with probability $p$ independently. Incorporating various strengthenings the phase transition can be summarised as follows. (All asymptotic statements are with respect to $n\to\infty$ and by \emph{whp} we abbreviate ``with probability $\to1$''.)
\begin{theorem}[Bollob\'as~\cite{Bollobas84}; \L uczak~\cite{Luczak90}]\label{thm:graph:giant}
Let $\eps=\eps(n)>0$ be a real function satisfying $\eps\to0$ and $\eps^3n\to\infty.$ 
\begin{enumerate}[(a)]
\item If $p=\frac{1-\eps}{n}$, then whp all components in $G(n,p)$ have size at most $O(\eps^{-2}\log (\eps^3n))$;
\item  If $p=\frac{1+\eps}{n}$, then whp the size of the largest component in $G(n,p)$ is $(1\pm o(1))2\eps n$, while all other components have size at most $O(\eps^{-2}\log (\eps^3n))$. 
\end{enumerate}
\end{theorem}
As we continue to add edges one by one, more and more components are consumed by the giant component and eventually the graph becomes connected. In fact, Bollob\'as and Thomason \cite{BollobasThomason85} showed that this happens precisely at the moment when the last isolated vertex disappears -- thereby relating a \emph{global graph property} to its minimal \emph{local} obstruction. Denote the \emph{hitting time} of connectivity  by $\conntime$, i.e.\ $\conntime$ is the minimal $M$ such that $G(n,M)$ is connected, and the hitting time for the disappearance of the last isolated vertex by $\isoltime$.

\begin{theorem}[Bollob\'as and Thomason~\cite{BollobasThomason85}]\label{thm:graph:hitting}
Whp in $\graphprocess$ we have $\conntime=\isoltime$.
\end{theorem} 

\section{Evolution of random hypergraphs -- Main results}

Given an integer $k\ge 2$ a \emph{$k$-uniform hypergraph $H$} consists of a set $V$ of vertices and a set $E$ of edges, where each edge contains precisely $k$ vertices. In particular, $2$-uniform hypergraphs are simply graphs. Given an integer $1\le j\le k-1$ we say that two $j$-sets (sets of $j$ distinct vertices) $J$ and $J'$ are \emph{$j$-connected} if there is a sequence of edges $e_1,\dots,e_m$ such that $J\subset e_1,$ $J'\subset e_m$ and $\left|e_i\cap e_{i+1}\right|\ge j$ for all $1\le i\le m-1.$ A \emph{$j$-component} is a maximal set of pairwise $j$-connected $j$-sets. The hypergraph $H$ is \emph{$j$-connected} if every two $j$-sets are $j$-connected. A $j$-set is called \emph{isolated} if it is not contained in any edge. Note that connectivity in graphs corresponds to the case $k=2$ and $j=1$.

We consider the $k$-uniform random hypergraph process $\process$: The vertex set is $[n]$ and initially there are no edges present. In each step of the process, we add an edge for a $k$-set chosen uniformly at random from all $k$-sets that do not already form an edge. Instead of analysing this process directly we usually consider the $k$-uniform binomial random hypergraph $\hknp$ with vertex set $[n]$, where every $k$-set is an edge with probability $p$ independently. This way there are no dependencies between different edges. It is well-known that both models are very similar and results can be easily transferred from one to the other using standard techniques (e.g.\ \cite{JansonLuczakRucinskiBook}).

For any $k\ge 2$ the case of vertex-connectivity ($j=1$) is well-studied and results analogous to Theorem~\ref{thm:graph:giant} were obtained in~\cite{BehrischCojaOghlanKang10,BollobasRiordan12c,KaronskiLuczak02,SchmidtShamir85}. Theorem~\ref{thm:graph:hitting} has also recently been extended for vertex-connectivity in $k$-uniform hypergraphs~\cite{Poole14}.

However, for high-order connectivity ($j>1$) not much was known until recently. This is due to the fact that vertex-connectivity can usually be studied with tools which are very similar to those used for graphs. By contrast, analysing high-order connectivity is often significantly more sophisticated and thus these methods are usually not sufficient. Recently Cooley, Kang, and Person~\cite{CooleyKangPerson14} showed that the sharp threshold for the emergence of the giant component in $\hknp$ is $$\giantthreshold:=\frac{1}{\left(\binom{k}{j}-1\right)\binom{n}{k-j}}.$$ We strengthen this result and provide the asymptotic size of the largest $j$-component in the \emph{weakly supercritical regime} of the phase transition. 

\begin{thm}\label{thm:main}
Let $k\ge 2$ and $1 \le j \le k-1$ be integers. Let $\eps=\eps(n)>0$ be a real function satisfying $\eps\rightarrow 0$, $\eps^3 n^j\rightarrow \infty$ and $\eps^2 n^{1-2\delta} \rightarrow \infty$, for some constant $\delta >0$. Then whp the size of the largest $j$-component $\mathcal{L}_1$ in $\cH^k(n,(1+\eps)p_g)$ satisfies
 $$
 \left|\mathcal{L}_1\right|=(1\pm o(1))\frac{2\eps}{\binom{k}{j}-1}\binom{n}{j},
 $$ 
 while all other $j$-components contain $o(\eps n^j)$ $j$-sets.
\end{thm}

The proof of Theorem~\ref{thm:main} fundamentally uses a powerful tool, the \emph{smooth boundary lemma}, which provides insight into the structure of the unique largest $j$-component in the supercritical regime. The precise statement of this lemma requires a large amount of additional notation and thus we omit it and refer to~\cite{CooleyKangKoch15a} for the details. Instead we explain the notion of `smooth sets' on a more intuitive level in Section~\ref{sec:smooth}. 

An elegant application of the notion of smoothness arises when studying the threshold for $j$-connectivity in $\hknp.$ The key idea is that the giant component contains a smooth set of large size. 

\begin{lemma}\label{lem:smoothsubset}
Let $k\ge2$ and $1 \le j \le k-1$ be integers. Let $\gamma=\gamma(n)>0$ be a real function satisfying $\gamma\rightarrow 0$ and $\gamma^3 n\rightarrow \infty$. Then whp the unique largest $j$-component  $\mathcal{L}_1$ in $\cH^k(n,(1+\gamma)p_g)$ contains a subset $S\subset\mathcal{L}_1$ of at least $\gamma^3 n^{j}$ $j$-sets with the following property for all integers $0\le \ell<j$:
\begin{center}
Every $\ell$-set $L\subset [n]$ is contained in $(1\pm o(1))\frac{\left|S\right|}{\binom{n}{j}}\binom{n}{j-\ell}$ $j$-sets of $S$.
\end{center}
\end{lemma}

Based on Lemma~\ref{lem:smoothsubset} we provide an elementary proof showing that the hypergraph process $\process$ becomes $j$-connected exactly at the moment when the last isolated $j$-set disappears. Let $\conntime=\conntime(k,j)$ be the hitting time for $j$-connectivity in the random hypergraph process $\process$ and let $\isoltime=\isoltime(k,j)$ denote the time-step of $\process$ in which the last isolated $j$-set disappears. 

\begin{thm}\label{thm:hittingtime}
Let $k\ge 3$ and $1\le j \le k-1$ be integers. Then whp in  $\process$ we have $\conntime = \isoltime$.
\end{thm}

It follows that, in the binomial random hypergraph $\hknp$, the properties of being $j$-connected and having no isolated $j$-sets share the common sharp threshold $$\connthreshold:=\frac{j\log n}{\binom{n}{k-j}}.$$ In fact we obtain a slightly stronger result.

\begin{thm}\label{thm:connthres}
Let $k\ge 3$ and $1\le j \le k-1$ be integers. Let $\omega=\omega(n)>0$ be a real function satisfying $\omega\to\infty$ and $\omega=o(\log n)$.
\begin{itemize}
\item If $p=\frac{j\log n + \omega}{\binom{n}{k-j}}$, then whp $\hknp$ contains isolated $j$-sets (and is therefore not $j$-connected);
\item If $p=\frac{j\log n - \omega}{\binom{n}{k-j}}$, then whp $\hknp$ is $j$-connected (and therefore contains no isolated $j$-sets).
\end{itemize}
\end{thm} 

\section{Proof outlines}

\subsection{Smooth sets in large components}\label{sec:smooth}

The intuition behind the smooth boundary lemma is the following: For an arbitrary $j$-set $J$, we explore its component via a breadth-first search process, i.e.\ generation by generation. If the component of $J$ happens to be large, the sizes of the generations have a tendency to grow and thus most generations should have a `reasonable' size already early on in the process. However, once the generations are not too small, random fluctuations should start to even themselves out over time. Thus generations should begin to look `smooth' in the sense that any set $L\subset[n]$ of at most $j-1$ vertices is contained in approximately the `right' number of $j$-sets of any smooth generation. Usually we will only run the exploration process for as long as necessary, resulting in a \emph{partial component}. If so, the last generation that is discovered before the process was stopped is called the \emph{boundary} and will be of special interest in the proof of Theorem~\ref{thm:main}, since it contains all the $j$-sets that are still \emph{active}. This will be discussed further in Section~\ref{sec:giant}. 

In other scenarios, the boundary does not play such a crucial role or may not be large enough on its own, in which case we may consider the smooth set given by the union of  all smooth generations in the partial component. It turns out that this union contains almost all $j$-sets of the partial component. As an immediate consequence we obtain Lemma~\ref{lem:smoothsubset}. The details can be found in~\cite{CooleyKangKoch15a,CooleyKangKoch15b}.

\subsection{Emergence of the giant component}\label{sec:giant}

First we study the number of $j$-sets in `large' components. For a given $j$-set we explore its component via a breadth-first search process and approximate this process by a supercritical branching process. In order to control the second moment we have to study two exploration processes and make sure that one of them being large does not increase the probability of the other becoming large too much. We first run one exploration process, stop it as soon as we know whp that it will grow large, and then consider a subprocess of the second exploration process in which no $k$-sets containing a $j$-set from the first (partial) component is present. The smooth boundary lemma ensures that the resulting process is still close to a (sufficiently) supercritical branching process. Therefore the number of $j$-sets in large components is concentrated around its expectation and Theorem~\ref{thm:main} follows by a sprinkling argument. The details can be found in~\cite{CooleyKangKoch15a}.

\subsection{Hitting times and connectivity threshold}\label{sec:hitting}

It is convenient to split the proof of Theorem~\ref{thm:hittingtime} into two parts: First we analyse the structure of the binomial random hypergraph $\hknp$ in the regime where the expected number of isolated $j$-sets begins to vanish. In particular, we show that whp in this regime the random hypergraph consists only of one non-trivial component and isolated $j$-sets. Then we transfer this structure to the corresponding time-range in the random hypergraph process $\process$ by classical contiguity results.

For integer-valued random variables $X_1,X_2,\dots$ and $Y$, we say \emph{$X_n$ converges in distribution to $Y$}, denoted by $X_n \stackrel{d}{\longrightarrow} Y$, if  we have $\Pr(X_n=i)\rightarrow \Pr(Y=i)$ for every integer $i$. We apply the Chen-Stein method for Poisson-approximation to the number $D_s$ of $j$-sets of degree $s\ge 0$ (i.e.\ $j$-sets that are contained in precisely $s$ edges). 
\begin{theorem}\label{thm:poisson}
Suppose $p=\frac{j\log n+s\log\log n +c_n}{\binom{n}{k-j}}$ for some real function $c_n=o(\log n)$. Then we have, for any integer $s\ge 0,$ 
\begin{center}
\begin{tabular}{rll}
$(i)$ & $D_s=0$ whp & if $c_n\rightarrow \infty$;\\
$(ii)$ & $D_s \stackrel{d}{\longrightarrow} \Po\left(\frac{j^s e^{-c}}{j!s!}\right)$ & if $c_n\rightarrow c\in\Real$;\\
$(iii)$ & $D_s\to \infty$ whp & if $c_n \to -\infty$.
\end{tabular}\end{center}
\end{theorem}
Set $p_1:=\frac{j\log n -\log\log n}{\binom{n}{k-j}}$  and $p_2:=\frac{j\log n +\log\log n}{\binom{n}{k-j}}$ and note that by Theorem~\ref{thm:poisson} with $s=0$ we know that whp $\cH^k(n,p_1)$ contains isolated $j$-sets, while $\cH^k(n,p_2)$ whp does not. Next set $p_0:=(1+\frac{1}{\log\log\log n})\giantthreshold$ and note that by Lemma~\ref{lem:smoothsubset} whp the unique largest component $\mathcal{L}_1$ in $\cH^k(n,p_0)$ already contains a smooth subset $S$ which contains every $(j-1)$-set at least $n/(2\log\log\log n)^3$ times. 

For any $p\in[p_1,p_2]$ we may construct $\hknp$ in a \emph{two-round exposure} such that $\hknp=\cH^k(n,p_0)\cup\cH^k(n,p^*)$, where $p^*=\frac{p-p_0}{1-p_0}=(1+o(1))p$, since $p_0\ll p$. In particular this means that we may assume that $\mathcal{L}_1$ and the smooth set $S$ are contained in $\hknp$. A first moment argument shows that whp $\hknp$ does not contain any non-trivial components consisting of at most $\log\log n$ $j$-sets. However, for much larger components this argument would not be sufficient, since the number of non-isomorphic hypergraphs consisting of a given number of $j$-sets grows too rapidly. This is where the smooth set $S$ comes into play. It allows us to give an upper bound for the expected number of large components that are not connected to $\mathcal{L}_1$ in $\cH^k(n,p^*)$ by calculating the expected number $X$ of \emph{connected subsets} of size approximately $\log\log n$, which are not connected to $S$ in $\cH^k(n,p^*)$. Since $S$ is smooth the number of forbidden edges can be bounded \emph{uniformly}. It follows that $X=o(1)$ and hence, by Markov's inequality, whp only the component containing $S$ is non-trivial. It remains to observe that the error-bounds in all these statements are strong enough to transfer them to $\process$ and apply a union bound over $2\log\log n \frac{\binom{n}{k}}{\binom{n}{k-j}}=\Theta(n^j\log\log n)$ time-steps.

Theorem~\ref{thm:connthres} is an immediate corollary of Theorem~\ref{thm:hittingtime} and Theorem~\ref{thm:poisson}. The details can be found in~\cite{CooleyKangKoch15b}.

\printbibliography

\end{document}